\documentstyle[amscd,amssymb,verbatim,11pt]{amsart}

\theoremstyle{plain}
\newtheorem{Prop}{Proposition}[section]
\newtheorem{Thm}[Prop]{Theorem}
\newtheorem{Cor}[Prop]{Corollary}

\newtheorem{Lem}[Prop]{Lemma}
\newtheorem{notation}[Prop]{Notation}

\theoremstyle{definition}
\newtheorem{Def}[Prop]{Definition}
\newtheorem{Idea}[Prop]{Idea}

\theoremstyle{remark}
\newtheorem{Rem}[Prop]{Remark}
\newtheorem{Problem}[Prop]{\bf Problem}

\def\dim{\mathop{\roman{dim}}}
\def\int{\mathop{\roman{int}}}

\def\1{^{-1}}

\def\Z{{\bold Z}}
\def\Q{{\bold Q}}

\def\p{{\bold p}}
\def\q{{\bold q}}
\def\B{\cal B_{\cal G}}

\def\dim{\text{dim}}

\def\Tor{\text{Tor}}

\def\ExD{\text{ext--dim}}

\def\int{\text{Int}}

\numberwithin{equation}{section}
\begin{document}

\title[Extension theory of infinite symmetric products]%
   {
Extension theory of infinite symmetric products}

\author{Jerzy Dydak}
\date{April 17, 2004}


\thanks{ Research supported in part by a grant
 DMS-0072356 from the National Science Foundation.}

\begin{abstract}

We present an approach to cohomological dimension theory
based on infinite symmetric products and on the
 general theory of dimension called
the extension dimension. The notion of the extension dimension $\ExD(X)$
was introduced
by A.N.Dranishnikov \cite {D$_5$} in the context
of compact spaces and CW complexes. 
 This paper investigates extension types of
infinite symmetric products $SP(L)$.
One of the main ideas of the paper is to treat
$\ExD(X)\leq SP(L)$ as the fundamental concept of
cohomological dimension theory
instead of $\dim_G(X)\leq n$. In a subsequent paper \cite{Dy$_6$}
we show how properties of infinite symmetric products
 lead naturally
to a calculus of graded groups which implies most of classical results of
the cohomological dimension. The basic notion in \cite{Dy$_6$} is that of homological
dimension of a graded group which allows for simultanous treatment
of cohomological dimension of compacta and extension properties
of CW complexes.
\par
We introduce cohomology
 of $X$ with respect to $L$ (defined as homotopy groups of the function space
$SP(L)^X$).
 As an application of our results we characterize all countable
groups $G$ so that the Moore space $M(G,n)$ is of the same
extension type as the Eilenberg-MacLane space $K(G,n)$.
Another application is characterization of infinite symmetric
products of the same extension type as a compact (or finite-dimensional
and countable) CW complex.

\end{abstract}

\maketitle

\medskip
\medskip
\tableofcontents

\section{Introduction}

\begin{notation} \label{XXXi.0} Throughout the paper $K$, $L$, and $M$ are
reserved for CW complexes. $X$ and $Y$ are general topological spaces
(quite often compact or compact metrizable).
We will frequently omit
coefficients in the case of
integral homology and cohomology. Thus,
 $H_n(K;\Z)$ will be shortened to $H_n(K)$ and $H^n(X;\Z)$ will be
shortened to $H^n(X)$.
\end{notation} 

Recall that $K$ is
an {\bf absolute extensor} of $X$ (this is denoted by $K\in AE(X)$)
if any map $f:A\to K$, $A$ closed in $X$,
extends over $X$.
That concept creates a partial order
on the class of CW complexes. Namely, $K\leq L$ if $K\in AE(X)$
implies $L\in AE(X)$ for every compact space $X$.
The partial order induces an equivalence relation on the class of all CW
complexes. The equivalence class $[K]$ of $K$ is called its
{\bf extension type}. 
\begin{Def} \label{XXXi.ext} A CW complex $K$ is called
the {\bf extension dimension} of a compact space $X$
(notation: $K=\ExD(X)$) if $K$ is a minimum
of the class $\{L\mid L\in AE(X)\}$.
\end{Def} 

\begin{Thm}[Dranishnikov Duality Theorem \cite {D$_5$}] \label{XXXi.1} 
Extension
dimension of compact spaces exists and each extension type
is equal to the extension dimension of some compact space.
\end{Thm} 

The concept of extension dimension generalizes both
the covering dimension $\dim(X)$ and the cohomological
dimension $\dim_G(X)$ with respect to an Abelian group $G$.
Indeed, $\dim(X)\leq n$ is equivalent to $\ExD(X)\leq S^n$
and $\dim_G(X)\leq n$ is equivalent to $\ExD(X)\leq K(G,n)$.

\par The theory of extension dimension is mostly
geometric in nature (see Section 2).
We introduce algebra to it following the basic idea of
\cite{Git}, where algebraic topology is outlined
via properties of infinite symmetric products $SP(K)$.
Thus, in this paper we show that the relation $SP(K)\leq SP(L)$
is of purely algebraic nature. We analyze it by generalizing
the connectivity index of Shchepin \cite{S} to the concept
of homological dimension of CW complexes.
To analyze relation $\ExD(X)\leq SP(L)$ we introduce
the concept of cohomology groups $H^\ast(X;L)$
of $X$ with coefficients in a CW complex $L$ (see Section 4).
Those cohomology groups have natural formulae facilitating
proofs and applications.
We show in Section 6 that the class $\{SP(L)\mid \ExD(X)\leq SP(L)\}$
has a minimum which
should be interpreted as the coefficient-free cohomological dimension
of $X$.
\par In Section 8 we dualize the connectivity index and use it
to derive algebraic implications of $\ExD(X)=\ExD(Y)$.
\par In a subsequent paper \cite{Dy$_6$} we will explain
that Bockstein Theory plays the role of homological algebra
in algebraic topology. In this paper we use Bockstein Theory
to give necessary and sufficient conditions for $SP(L)$
to have the same extension type as an Eilenberg-MacLane
space $K(G,n)$ (see Section 7).
Later on (Section 9) we characterize extension types
$[SP(L)]$ containing compact (respectively, countable and
finite-dimensional) CW complexes. That characterization
generalizes all the previously known results about
different extension types.

The last part of the paper (Section 9) is devoted to
comparison of extensional properties of $M(G,n)$ and $K(G,n)$.

\begin{Def} \label{XXXi.5} Suppose $G$ is an Abelian group and $n\ge 1$ is
an integer.
By $M(G,n)$ we will denote {\bf the Moore space}, i.e. a CW complex
$K$ so that $\tilde H_n(K;\Z)=G$ and $\tilde H_m(K;\Z)=0$ for $m\ne n$.

\end{Def}

More precisely, Moore spaces discussed in this paper are constructed as
follows.
Choose a short exact sequence $0\to F_1\to F\to G\to 0$ so that $F$ is free
Abelian.
Let $L$ be the wedge of $n$-spheres enumerated by generators of $F$. Attach
$(n+1)$-cells to $L$ enumerated by generators of $F_1$ via characteristic maps
corresponding to $F_1\to F$. In particular, such Moore spaces are
finite-dimensional
and one has a map $M(\pi_n(K),n))\to K$ 
(provided $\pi_1(K)$ is Abelian if $n=1$) inducing isomorphism of $n$-th
homotopy groups
for any space $K$.

\par

We will use the following generalization
of the relation $K\leq L$ (see \cite{D-D$_1$}):

\begin{Def} \label{XXXi.2} Suppose $K$ and $L$ are CW complexes and $\cal C$
is a class of spaces. $K\leq_{\cal C}L$ means that if $X\in\cal C$ and
$K\in AE(X)$,
then $L\in AE(X)$. In particular $K\leq_X L$ means $K\leq_{\{X\}} L$.
\end{Def} 

The author would like to thank Akira Koyama for discussions on the subject
of the paper, and E.Shchepin for giving a series of excellent talks
(in Russian) during the workshop \lq Algebraic ideas in dimension
theory\rq \ held in Warsaw (Fall 1998). Those talks were the starting
point in author's understanding of paper \cite{S}
which eventually led to the ideas developed in this paper.
\par Some of the results of this paper were circulated
in unpublished notes \cite{Dy$_7$}.

\section{Geometry of extension theory }

The purpose of this section is to list results of extension
theory which involve no algebraic computations.


\begin{Prop} \label{XXX2.1} (see \cite{D-D$_1$},\cite{D$_1$}) Suppose $K$
is a CW complex and $X$ is metrizable. If $K\in AE(X)$, then every map
$f:X\to \Sigma(K)$ from $X$ to the suspension of $K$ is null-homotopic.

\end{Prop}


\begin{Thm} \label{XXX2.2} (see \cite{D$_4$}) 
If $X=\bigcup\limits_{n=1}^\infty X_n$, $X_n$ is a closed subset of $X$ and $K\in
AE(X_n)$
for all $n$, then
$K\in AE(X)$ provided $X$ is normal and $K\in ANE(X)$.

\end{Thm}


\begin{Thm} \label{XXX2.3} (see \cite{Wa} and \cite{S}) 
Suppose $K$ is a CW complex. If $K\in AE(X)$, then
$K\in AE(Y)$ for every $Y\subset X$ if $X$ is metrizable.

\end{Thm}


\begin{Thm} \label{XXX2.4} (see \cite{Dy$_3$}) 
If $X=A\cup B$ is metrizable, and $K\in AE(A)$, $L\in AE(B)$
are CW complexes,
then $K\ast L\in AE(X)$.

\end{Thm}


\begin{Thm} \label{XXX2.5} (see \cite{D$_4$}) 
Suppose $K$ and $L$
are countable CW complexes.
If $K\ast L\in AE(X)$ and $X$ is a compactum,
then there is a $G_\delta$-subset $A$ of $X$ such that
  $K\in AE(A)$ and $L\in AE(X-A)$.

\end{Thm}


\begin{Thm} \label{XXX2.6} (see \cite{O$_1$}) 
Suppose $K$ is a countable CW complex. If $K\in AE(X)$ and $X$ is a subset
of a metric separable space $X'$,
then there is a $G_\delta$-subset $A$ of $X'$ containing $X$ such that
  $K\in AE(A)$.

\end{Thm}

\begin{Thm} \label{XXX2.11} (see \cite{D-D$_1$}) Suppose $X$ is compact or 
metrizable and $K$ is
a pointed connected CW complex. The following conditions are
equivalent:
\par{1.} $K\in AE(X\times I)$.
\par{2.} $\Omega(K)\in AE(X)$, where $\Omega(K)$ is the loop-space
of $K$.
\par{3.} $\lbrack X/A,K\rbrack=0$ for all closed subsets $A\ne \emptyset$ of $X$.
\par{4.} $K\in AE(\Sigma(X))$.

\end{Thm} 

\begin{Thm} \label{XXX2.9} (see \cite{D$_3$}) 
Suppose $K$ is a CW complex. If $X$ is finite-dimensional and
$\prod\limits_{i=1}^{\infty}K(\pi_i(K),i)\in AE(X)$, then $K\in AE(X)$.

\end{Thm}

\section{Transition to algebra in extension theory}

 Given a space
$X$ and $k>0$, the $k$-th {\bf symmetric product}
$SP^k(X)$ of $X$ is the  space of orbits of the action of
the symmetric group $S_k$ on $X^k$.  Points of  $SP^k(X)$ can be written
in the form $\sum\limits_{i=1}^{k} x_i$. $SP^k(X)$ is equipped
 with the quotient topology given by the natural  map $\pi :X^k\to SP^k(X)$.
If $X$ is
metrizable, then $\pi :X^k\to SP^k(X)$ is both open and closed
(see p. 255 of \cite {D-T}), so $SP^k(X)$ is metrizable, too 
(use 4.4.18 of \cite {En}).
\par\noindent If $X$ has a base point $a$, then $SP^k(X)$
has $\sum\limits_{i=1}^{k} a$ as its base point. Notice that there is
a natural inclusion $i:SP^n(X)\to SP^k(X)$ for all $n<k$.
It is given by the formula $i(\sum\limits_{i=1}^{n}
x_i)=\sum\limits_{i=1}^{n} x_i+(k-n)a$.
In this way, points of the form  $\sum\limits_{i=1}^{n} x_i$, $n<k$, can
be considered as belonging to $SP^k(X)$.
\par\noindent The direct limit of $SP^2(X)\to \ldots \to SP^n(X)\to \ldots$
is denoted by
$SP(X)$ and is called {\bf the infinite symmetric product}
(see \cite {D-T} or \cite{Git}, p.168).
\par

The main property of the infinite symmetric product is expressed
in the following result.


\begin{Thm}[Dold-Thom Theorem \cite {D-T}] \label{XXX3.1} If $K$ is a pointed CW complex, then
the natural inclusion $i:K\to SP(K)$
induces an isomorphism $\tilde H_i(K)\to \pi_i(SP(K))$.

\end{Thm} 

 The meaning of Dold-Thom Theorem is that one can define
singular homology groups geometrically, without the apparatus of
homological algebra,
as homotopy groups of infinite symmetric products (see \cite{Git}).

Of major importance to us is the following result of Dranishnikov:
\begin{Thm} \label{XXX2.10a} (see \cite{D$_5$}) 
Suppose $K$ is a CW complex. If $X$ is compact and $K\in AE(X)$, then
$SP(K)\in AE(X)$. 
\end{Thm}

Since $SP(K)$ is homotopy equivalent
to the weak product of Eilenberg-MacLane spaces
$K(\tilde H_i(K),i)$ (see \cite{Git}, Corollary 6.4.17 on p.223)
one has the following:
\begin{Thm} \label{XXX2.10b} (see \cite{D$_5$}) 
Suppose $K$ is a CW complex. If $X$ is compact, then
$SP(K)\in AE(X)$
is equivalent to $K(\tilde H_i(K),i)\in AE(X)$ for all $i\ge 0$. 
\end{Thm}

The following result of Dranishnikov uses a high level
of algebraic arguments and can be viewed as an analog
of the Hurewicz Theorem in Extension Theory.

\begin{Thm} \label{XXX2.10} (see \cite{D$_5$}) 
Suppose $K$ is a CW complex. If $X$ is compact
finite-dimensional, $K$ is simply connected, and
$SP(K)\in AE(X)$, then $K\in AE(X)$.
\end{Thm}

\par Let us move to algebraic concepts associated with extension theory
by employing the connectivity index introduced by
 E.Shchepin \cite{S}. 

\begin{Def} \label{XXX1.3} (\cite{S},p.985) 
Suppose $K$ is a CW complex. Its {\bf connectivity index}
$cin(K)$ is either $\infty$ or a non-negative integer. $cin(K)=\infty$ means
that all reduced integral homology groups of $K$ are trivial.
$cin(K)=n$ means that $\tilde H_n(K;\Z)\ne 0$
and $\tilde H_k(K;\Z)=0$ for all $k<n$.
\end{Def} 

\begin{Prop} \label{XXX2.7} Suppose $K$ is a simply-connected
CW complex. The following numbers are equal:
\par\par{1.} $cin(K)$, 
\par\par{2.} the supremum of all $n\ge 0$ so that
$K\in AE(S^n)$,
\par\par{3.} the supremum of all $n\ge 0$ so that
$K$ is homotopy $k$-connected for all $k<n$.
\end{Prop}

\begin{pf}  If $K\in AE(S^n)$, then any map $S^j\to K$, $j<n$,
extends over $S^n$ and must be null-homotopic. Hence
$H_j(K)=0$ for $j<n$.
If $H_j(K)=0$ for $j<n$, then $\pi_j(K)=0$ for $j<n$
and any map $f:A\to K$, $A$ closed in $S^n$, can be extended
over $S^n$ as follows: first extend it over a polyhedral
neighborhood $N$ of $A$, then keep extending over simplices
which are not contained in $N$ using induction over
the dimension of simplices.
 \end{pf}

\begin{Prop} \label{XXX2.8} Suppose $K$ and $L$ are pointed CW
complexes.
If $K$ is countable and $K\le_{X} L$ for all finite-dimensional compacta $X$,
then $cin(K\wedge M)\le cin(L\wedge M)$
for every CW complex $M$.
\end{Prop}

\begin{pf}  
\par Case 1: $M$ is countable.
\par 
Since suspending a CW complex pushes its connectivity index up by 1,
it suffices to prove \ref{XXX2.8} for $M=\Sigma^2 M'$, in which case $K\wedge M$
and $L\wedge M$ are simply connected and we may use \ref{XXX2.7}.
Notice that $cin(K\wedge M)=cin(K\ast M)-1$, where $K\ast M$
is the join of $K$ and $M$.
Suppose $K\ast M\in AE(S^m)$. Express $S^m$ as $A\cup B$
so that $K\in AE(A)$, $M\in AE(B)$, and $A$ is $F_\sigma$ (see  \ref{XXX2.5}).
Now, $L\in AE(A)$ (see \ref{XXX2.2}) which implies 
$L\ast M\in AE(S^m)$ by \ref{XXX2.4}.
\par Case 2:  $M$ is arbitrary.
\par Suppose $m=cin(K\wedge M)> cin(L\wedge M)=n$.
There is an integral cycle $c$
in $H_n(L\wedge M)-\{0\}$. That cycle lies in $L\wedge M_1$ for
some finite subcomplex $M_1$ of $M$. 
Given countable subcomplex $M_i$ of $M$ construct countable subcomplex 
 $M_{i+1}\supset M_i$ of $M$ so that $H_k(K\wedge M_i)\to H_k(K\wedge M_{i+1})$
is trivial for $k<m$.
Let $M'$ be the union of all $M_i$.
By Case 1, $m\leq cin(K\wedge M')\leq cin(L\wedge M')$
which means that $c$ becomes $0$ in $H_n(L\wedge M')$,
a contradiction.
 \end{pf}

\section{Cohomology with coefficients in a complex}

\par With the hindsight of homological algebra
one has the pairing
$$H_*(K)\wedge H_*(L)\to \pi_*(SP(K\wedge L))$$
for any two pointed CW complexes $K$ and $L$.
It corresponds to the well-known K\" unneth formula for homology
and we will give it a slightly non-traditional form.


\begin{Thm} \label{XXX3.2} If $K$ and $L$ are pointed
CW complexes, then
$$H_n(K\wedge L)\equiv \bigoplus\limits_{i+j=n}H_i(K;H_j(L)).$$

\end{Thm}

\begin{Rem} \label{XXX3.2.5}
 Theorem \ref{XXX3.2} is a direct consequence of the Universal
Coefficient Theorem (see \cite{Sp}, p.222) and the classical
K\" unneth Theorem (see \cite{Sp}, p.235). 
\end{Rem}

Notice that
$$\bigoplus\limits_{i+j=n}H_i(K;H_j(L))\equiv
\bigoplus\limits_{i+j=n}H_i(L;H_j(K))$$
as both groups are isomorphic to $H_n(K\wedge L)$. A natural question to
ask is if
$$\bigoplus\limits_{i}H^i(K;H_{n-i}(L))$$
or
$$\bigoplus\limits_{i}H_i(L;H^{n-i}(K))$$
have similar geometrical interpretation, and if they are isomorphic.
Knowing that the smash product
$K\wedge L$ is adjoint to the function space functor
one can speculate the following:


\begin{Idea} \label{XXX3.3}
There ought to be a dual pairing
$$H^*(X)\wedge H_*(L)\equiv \pi_*(SP(L)^X).$$

\end{Idea}

The above pairing should correspond to K\" unneth formula for cohomology.
It turns out such pairing exists for pointed compact spaces.
\par Our first idea is to introduce cohomology of a pointed space
$X$ via homotopy groups of function spaces $SP(K)^X$.


\begin{Def} \label{XXX3.4} Suppose $K$ is a pointed CW complex and $X$
is a pointed space. The {\bf cohomology} $H^n(X;K)$ of $X$ {\bf
with coefficients} in $K$ is defined as follows:
\par{1.} $H^n(X;K)=[X,SP(\Sigma^{n}(K)]$ if $n\ge 0$,
\par{2.} $H^n(X;K)=[\Sigma^{-n}X,SP(K)]$ if $n\leq 0$.
\end{Def}


\begin{Prop} \label{XXX3.5} If $K=M(G,n)$ is a Moore space, then $H^k(X;K)=H^{n+k}(X;G)$.

\end{Prop}

\begin{pf}  Notice that $\Sigma^r(K)=M(G,n+r)$ and
$SP(M(G,r))=K(G,r)$ for $r\ge 1$ (see 3.1). As $H^m(X;G)=[X,K(G,m)]$
one gets \ref{XXX3.5} immediately from the definition of $H^m(X;K)$.
 \end{pf}

It is shown in \cite{D-T} that $SP(K)$
is homotopy equivalent to $\bigcup\limits_{m=1}^\infty\prod\limits_{n=1}^m
K(H_n(K),n)$ for every connected CW complex $K$
(see also \cite{Git}, Corollary 6.4.17 on p.223). We will need a more general
result.


\begin{Prop} \label{XXX3.6} Suppose $K$ is a pointed CW complex
and $X$ is a pointed $k$-space. There is a weak homotopy equivalence
$i:\bigcup\limits_{m=1}^\infty\prod\limits_{n=1}^m
K(\pi_n(SP(K)^X),n)\to SP(K)^X$.

\end{Prop}

\begin{pf}  Notice any map $f:L\to SP(K)^X$ extends
to $F:SP(L)\to SP(K)^X$ if $L$ is a CW complex.
Indeed, one can define $F$ as follows:
$F(\sum\limits_{i=1}^na_i)(x)=\sum\limits_{i=1}^nf(a_i)(x)$
for $\sum\limits_{i=1}^na_i\in SP(L)$ and $x\in X$.
Let $L$ be the wedge of $M(\pi_n(SP(K)^X),n)$, $n\ge 1$.
There is a map $f:L\to SP(K)^X$ so that
$\pi_n(f|M(\pi_n(SP(K)^X),n))$ is an isomorphism for each $n$
(see the discussion after  \ref{XXXi.5}). In particular, $H_k(L)\to \pi_k(SP(K)^X)$
is an isomorphism for all $k\ge 1$.
Pick an extension $F:SP(L)\to SP(K)^X$ of $f$.
Since $H_k(L)\to \pi_k(SP(L))$
is an isomorphism for all $k\ge 1$ and $H_k(L)\to \pi_k(SP(K)^X)$
is an isomorphism for all $k\ge 1$, $F$ is a weak homotopy
equivalence. As shown in \cite{D-T},
$SP(L)$
is homotopy equivalent to $\bigcup\limits_{m=1}^\infty\prod\limits_{n=1}^m
K(H_n(L),n)$ which is exactly the space
 $\bigcup\limits_{m=1}^\infty\prod\limits_{n=1}^m
K(\pi_n(SP(K)^X),n)$.
 \end{pf}


\begin{Cor} \label{XXX3.7} Suppose $K$ is a pointed connected
CW complex and $X$ is a pointed $k$-space. $SP(K)$ and
the spaces $\Omega^r(SP(\Sigma^r(K)))$
are homotopy equivalent for $r\ge 1$. In particular,
 $H^n(X,\Sigma^r(K))\equiv H^{n+r}(X;K)$
for all $n\in\Z$ and all $r\ge 1$.

\end{Cor}

\begin{pf}  Notice that $\pi_k(\Omega^r(SP(\Sigma^r(K))))=
\pi_{k+r}(SP(\Sigma^r(K)))=H_{k+r}(\Sigma^r(K))=H_k(K)$ for 
each $k\ge 0$.
By \ref{XXX3.6} (applied to $X=S^r$) there is is a weak homotopy equivalence
 $\bigcup\limits_{m=1}^\infty\prod\limits_{n=1}^m
K(H_{n}(K),n)\to \Omega^r(SP(\Sigma^r(K)))$. This map
is a homotopy equivalence as both spaces are homotopy
equivalent to CW complexes.
As in \cite{D-T}, 
$\bigcup\limits_{m=1}^\infty\prod\limits_{n=1}^m
K(H_{n}(K),n)$ is homotopically equivalent to $SP(K)$.
\par If $n\ge 0$, then
\par\noindent
$H^n(X;\Sigma^r(K))=[X,SP(\Sigma^{n+r}(K)]=H^{n+r}(X;K)$.
\par
If $n<0$, then $H^n(X;\Sigma^r(K))=[\Sigma^{-n}(X),SP(\Sigma^r(K))]$.
If $n+r\ge 0$, then
\par\noindent
$[\Sigma^{-n}(X),SP(\Sigma^r(K))]=
[X,\Omega^{-n}(SP(\Sigma^r(K)))]=[X,(SP(\Sigma^{n+r}(K)))]=
H^{n+r}(X;K)$.
If $n+r<0$, then
$[\Sigma^{-n}(X),SP(\Sigma^r(K))]=
[\Sigma^{-n-r}(X),\Omega^{r}(SP(\Sigma^r(K)))]$
\par\noindent
$=[\Sigma^{-n-r}(X),SP(K)]=
H^{n+r}(X;K)$.
 \end{pf}

The purpose of the next result is to generalize the well known
theorem of Cohen \cite{C}.


\begin{Thm} \label{XXX3.8} Suppose $K$ is a pointed CW complex
and $X$ is compact or metrizable. The following conditions are equivalent:
\par{1.} $SP(K)\in AE(X)$.
\par{2.} $H^n(X/A;K)=0$ for all $n>0$ and
all closed subsets $A$ of $X$.
\par{3.} $H^{1}(X/A;K)=0$ for
all closed subsets $A$ of $X$.

\end{Thm}

\begin{pf}  $H^n(X/A;K)$ was defined as $[X/A,SP(\Sigma^n(K))]$
for $n\ge 1$ (see \ref{XXX3.4}). \ref{XXX2.11} says that
$[X/A,SP(\Sigma(K))]=0$ for all closed subsets $A$ of $X$
if and only if $\Omega(SP(\Sigma(K)))\in AE(X)$. By \ref{XXX3.7},
$\Omega(SP(\Sigma(K)))$ is homotopy equivalent
to $SP(K)$ which proves 1)$\iff$ 3).
\par Clearly, 2) is stronger than 3).
\par Suppose 3) and 1) hold. If $n>1$, then (see \ref{XXX2.11})
$SP(\Sigma^n(K))\in AE(X)$
and $H^1(X/A;\Sigma^{n-1}(K))=0$ for all closed subsets $A$ of $X$
(as 1) is equivalent to 3) for all CW complexes $K$).
Now, \ref{XXX3.7} says $H^n(X/A;K)=H^1(X/A;\Sigma^{n-1}(K))$ which proves that 2) holds.
 \end{pf}

Notice that in \ref{XXX3.4} the homotopy groups of $(SP(K))^X$ correspond
to negative cohomology groups of $X$ with coefficients in $K$. It seems
unnatural
but we chose it that way in order to adhere to common practices.
However, it is time to break with tradition and achieve a better theory
that way.


\begin{Def} \label{XXX3.9} Suppose $h^*$ is a cohomology theory.
By ${}^rh^*$ we will
denote the {\bf reversed cohomology} defined via $^rh^m(X)=h^{-m}(X)$.

\end{Def}

The following result shows that using reversed cohomology one can
achieve similarity between homology and cohomology (compare \ref{XXX3.10} with 
 \ref{XXX3.2}).


\begin{Thm} \label{XXX3.10} If $K$ is a pointed CW complex
and $X$ is a pointed  compact space, then
$${}^rH^n(X;K)\equiv \bigoplus\limits_{i}{}^rH^i(X;H_{n-i}(K))$$
and 
$${}^rH^n(X;K)\equiv \bigoplus\limits_{i}H_i(K;{}^rH^{n-i}(X)).$$

\end{Thm}

\begin{pf}  If $Y$ is a pointed compact space and $L$ is a pointed
CW complex, then $[Y,SP(L)]$ is the direct
sum $\bigoplus\limits_{i}{}[Y,K(H_i(L),i)]$. Indeed,
$SP(L)$ is homotopically equivalent to 
$\bigcup\limits_{m=1}^\infty\prod\limits_{n=0}^m
K(H_{n}(L),n)$ (see \cite{D-T}) and any map from $Y$
(or any homotopy from $Y\times I$) to
$\bigcup\limits_{m=1}^\infty\prod\limits_{n=1}^m
K(H_{n}(L),n)$ is contained in $\prod\limits_{n=0}^m
K(H_{n}(L),n)$ for some $m$.
\par
If $n\ge 0$, then 
$${}^rH^{-n}(X;K)=[X,SP(\Sigma^{n}(K)]=$$
$$\bigoplus\limits_i
[X,K(H_{i-n}(K),i)]=\bigoplus\limits_{i}{}^rH^{-i}(X;H_{i-n}(K)).$$
\par
If $n< 0$, then ${}^rH^{-n}(X;K)=[\Sigma^{-n}X,SP(K)]$ is
 $\bigoplus\limits_i
[\Sigma^{-n}(X),K(H_{i}(K),i)]$ which is the same as
$\bigoplus\limits_{i}H^i(X;H_{i-n}(K))=\bigoplus\limits_{i}{}^rH^{-i}(X;H_{i-n}(
K)).$
\par By the Universal Coefficient Theorem for homology (see Theorem 14
in \cite{Sp} on p.226)
$\bigoplus\limits_{i}H_i(K;{}^rH^{n-i}(X))$ is isomorphic
to 
$$G_1=\bigoplus\limits_{i}(H_i(K)\otimes {}^rH^{n-i}(X)\oplus H_{i-1}(K)*{}^rH^{n-i}(X))$$
($G*G'$ is the torsion product of $G$ and $G'$).
By the Universal Coefficient Theorem for cohomology (see Statement 5
in \cite{K} on p.4)
$\bigoplus\limits_{i}{}^rH^i(X;H_{n-i}(K))$ is isomorphic
to 
$$G_2=\bigoplus\limits_{i}(H_{n-i}(K)\otimes {}^rH^{i}(X)\oplus H_{n-i}(K)*{}^rH^{i-1}(X)).$$
Notice that $G_1$ is isomorphic to $G_2$ (change $i$ to $n-i$ in the formula for $G_1$).
 \end{pf}

As a simple consequence of the fact that $K^Y$ is homotopy equivalent to a CW
complex if $K$ is a CW complex and $Y$ is compact, we get the following
version of K\" unneth
Formula.


\begin{Thm} \label{XXX3.11} Suppose $K$ is a pointed CW complex.
If $X$ 
and $Y$ are pointed compact spaces, then
$${}^rH^n(X\wedge Y;K)\equiv \bigoplus\limits_{i}{}^rH^i(X;{}^rH^{n-i}(Y;K)).$$

\end{Thm}

\begin{pf}  Suppose $n\ge 0$.
${}^rH^n(X\wedge Y;K)=[\Sigma^n(X\wedge Y),SP(K)]=
[\Sigma^n(X),SP(K)^Y]$. Since $SP(K)^Y$ is homotopy equivalent
to a CW complex, \ref{XXX3.6} says that it is homotopy equivalent to
$\bigcup\limits_{m=1}^\infty\prod\limits_{i=0}^m
K({}^rH^{i}(Y;K),i)$. Thus,
$${}^rH^n(X\wedge Y;K)=\bigoplus\limits_iH^i(\Sigma^n(X);{}^rH^{i}(Y;K))=$$

$$\bigoplus\limits_iH^{i-n}(X;{}^rH^{i}(Y;K))=\bigoplus\limits_i{}^rH^{n-i}(X;{}^r
H^{i}(Y;K)).$$
\par If $n<0$, then one reduces it to the case $n=0$ by observing
that ${}^rH^n(Z;K)={}^rH^0(Z;\Sigma^{-n}(K))$ for every pointed $k$-space
$Z$ (see \ref{XXX3.7}).
 \end{pf}


\begin{Cor} \label{XXX3.12} If $K$ is a pointed CW complex
and $X$ is a pointed compact space, then
the following conditions are equivalent ($m$ is an integer):
\par{1.}
$H^n(X;K)=0$ for all $n\ge m$.
\par{2.}
$H_i(K;H^n(X))=0$ for all $i\leq n-m$.
\par{3.}
$H^{i}(X;H_{n}(K))=0$ for all $i\ge n+m$.

\end{Cor}

\begin{pf}  Notice that $H^n(X;K)={}^rH^{-n}(X;K)$
is isomorphic to $\bigoplus\limits_{i}H_i(K;{}^rH^{-n-i}(X))$ (see \ref{XXX3.10})
which is the same as $\bigoplus\limits_{i}H_i(K;H^{n+i}(X))$.
That means $H^n(X;K)=0$ for all $n\ge m$ is equivalent to
$H_i(K;H^{n+i}(X))=0$ for all $n\ge m$ which is the same as saying that
$H_i(K;H^n(X))=0$ for all $i\leq n-m$.
Similarly, $H^n(X;K)={}^rH^{-n}(X;K)$
is isomorphic to $\bigoplus\limits_{i}{}^rH^i(X;H_{-n-i}(K))$ (see \ref{XXX3.10})
which is the same as $\bigoplus\limits_{i}H^{-i}(X;H_{-n-i}(K))$.
That means $H^n(X;K)=0$ for all $n\ge m$ is equivalent to
$H^{-i}(X;H_{-n-i}(K))=0$ for all $n\ge m$ which is the same as saying
$H^{i}(X;H_{n}(K))=0$ for all $i\ge n+m$.
 \end{pf}


\begin{Cor} \label{XXX3.13} Suppose $K$ is a pointed CW complex
and $X$ is a compact space. The following conditions are equivalent:
\par{1.} $SP(K)\in AE(X)$.
\par{2.} $H_i(K;H^n(X;A))=0$ for all $i<n$ and all closed subsets $A$ of $X$.
\par{3.} $H^i(X/A;H_n(K))=0$ for all $n<i$ and all closed subsets $A$ of $X$.

\end{Cor}

\begin{pf}  By \ref{XXX3.8}, $SP(K)\in AE(X)$ if and only if
$H^n(X/A;K)=0$ for all $n\ge 1$ and all closed subsets $A$ of $X$. 
Using \ref{XXX3.12} one gets that
$SP(K)\in AE(X)$ if and only if $H_i(K;H^n(X,A))=0$ for all $i\leq n-1$.
That proves 1)$\iff$ 2). 2)$\iff$ 3) follows from \ref{XXX3.12} applied to $m=1$.
 \end{pf}

\section{Homological dimension of CW complexes}

\ref{XXX3.13} suggests that $K\leq L$ should have
 some algebraic implications. The next result specifies the nature
of those implications.


\begin{Thm} \label{XXX3.14} Suppose $n>0$, $G$ is an Abelian group,
$K,L$ are CW complexes, and $\cal C$ is a class of spaces
containing all finite-dimensional
compacta. If $K$ is countable, $K\leq_{\cal C}L$, and $\tilde H_i(K;G)=0$ for all $i\leq n$, then
$\tilde H_i(L;G)=0$ for all $i\leq n$.

\end{Thm}

\begin{pf}  
Make $K$ and $L$ pointed CW complexes and switch from the reduced homology
to ordinary homology of pointed CW complexes.
By \ref{XXX3.2}, $H_i(A\wedge M(G,1))=H_{i-1}(A;G)$ for any pointed
CW complex $A$, where $M(G,1)$ is
the Moore space. Thus $cin(K\wedge M(G,1))\ge n+2$. By \ref{XXX2.8},
$cin(L\wedge M(G,1))\ge n+2$ which is the same as 
$\tilde H_i(L;G)=0$ for all $i\leq n$.
 \end{pf}

Theorem \ref{XXX3.14} suggests a new concept
of {\bf homological dimension} $\dim_G(K)$
of a pointed CW complex.

\begin{Def} \label{XXX3.15a} Suppose $K$ is a pointed CW complex and $G$
is an Abelian group. $\dim_G(K)=n < \infty$ means that $H_i(K;G)=0$ for all
$i < n$ and
$H_n(K;G)\ne 0$. $\dim_G(K)=\infty$ means that $H_i(K;G)=0$ for 
all $i$.
\end{Def} 

\begin{Rem} \label{XXX3.15b}
Notice that the above concept generalizes the concept of connectivity
index. Indeed $cin(K)=\dim_\Z(K)$ for all pointed CW complexes $K$.
\end{Rem}

Definition \ref{XXX3.15a} suggests a new partial order on the class of CW complexes:


\begin{Def} \label{XXX3.15} Suppose $K$ and $L$ are CW complexes and $G$
is an Abelian group. $K\leq_GL$ means $\dim_G(K)\leq \dim_G(L)$.
\par
If $\cal G$
is a class of Abelian groups, then $K\leq_{\cal G}L$ means that
$K\leq_GL$ for all $G\in\cal G$. $K\sim_{\cal G}L$ means that
$K\leq_{\cal G}L$ and $L\leq_{\cal G}K$.
\par $K\leq_{Gr} L$ means that
$K\leq_GL$ for all Abelian groups $G$. $K\sim_{Gr} L$ means that
$K\leq_{Gr} L$ and $L\leq_{Gr} K$.

\end{Def}


\begin{Cor} \label{XXX3.16} $K\sim_{Gr} SP(K)$ for each pointed
CW complex $K$.

\end{Cor}

\begin{pf}  
Notice that $H_n(K)$ is a direct summand of $H_n(SP(K))$ for each $n$.
 Indeed, $SP(K)$ is homotopy equivalent
to $\bigcup\limits_{m=1}^\infty\prod\limits_{n=1}^m
K(H_n(K),n)$ and each $K(H_n(K),n)$ (whose $n$-th homology group is 
$H_n(K)$) is a retract of $SP(K)$.
If $H_k(SP(K);G)=0$ for $k<n$, then it amounts to
$H_k(SP(K))\otimes G=H_{k-1}(SP(K))*G=0$ for $k<n$
(in view of the Universal Coefficient Theorem). Therefore
$H_k(K)\otimes G=H_{k-1}(K)*G=0$ for $k<n$ and $H_k(K;G)=0$ for $k<n$.
This proves $SP(K)\leq_{Gr} K$.
\par To prove
$K\leq_{Gr} SP(K)$ notice that, for $K$ countable, it follows 
from \ref{XXX3.14} and  \ref{XXX2.10a}.
Suppose $H_k(K;G)=0$ for $k\leq n$ and there is $c\in H_n(SP(K);G)-\{0\}$.
Find a finite subcomplex $K_1$ of $K$ and a countable subgroup $G_1$ of $G$
such that $c$ belongs to the image of $H_n(SP(K_1);G_1)\to H_n(SP(K);G)$.
By induction find countable subcomplexes $K_1\subset K_2\subset \ldots$
of $K$ and countable subgroups $G_1\subset G_2\subset\ldots$ of $G$
so that $H_k(K_i;G_i)\to H_k(K_{i+1};G_{i+1})$ is trivial for all $i$
and all $k\leq n$. Let $K'$ be the union of all $K_i$ and let $G'$ be the union of all $G_i$.
Notice that $H_k(K';G')=0$ for $k\leq n$. Therefore $H_n(SP(K');G')=0$
contradicting $c\ne 0$.
 \end{pf}


\begin{Thm} \label{XXX3.17} Suppose $K$ and $L$ are connected CW
complexes and $K$ is countable.
Consider the following conditions:
\par{1.} $K\leq_{Gr} L$.
\par{2.} $SP(K)\leq_X\ SP(L)$ for all compact $X$.
\par{3.} $cin(K\wedge M)\le cin(L\wedge M)$ for each complex $M$.
\par{4.} $K\leq_X\ L$ for all finite-dimensional compacta $X$.
\par Condition 4 implies Condition 1. Conditions 1-3 are equivalent. 
If $L$ is simply connected, then Conditions 1-4 are equivalent.

\end{Thm}

\begin{pf}  4)$\implies $1) follows from \ref{XXX3.14}.
\par
3) $\implies$ 1). Suppose $H_k(K;G)=0$ for $k<n$.
Use $M=M(G,1)$ and \ref{XXX3.2} to conclude that $H_k(K\wedge M)=0$
for $k\leq n$. Therefore $H_k(L\wedge M)=0$
for $k\leq n$ which means $H_k(L;G)=0$
for $k<n$.
\par
1) $\implies$ 3). Suppose  $H_k(K\wedge M)=0$
for $k< n$. That means (see \ref{XXX3.2}) $H_i(K;H_j(M))=0$ if $i+j<n$.
Since $K\leq_{Gr} L$, $H_i(L;H_j(M))=0$ if $i+j<n$,
i.e. $H_k(L\wedge M)=0$
for $k< n$.
\par 2)$\implies $1) follows from \ref{XXX3.14} and \ref{XXX3.16}.
\par  1)$\implies $2). Suppose $SP(K)\in AE(X)$.
\ref{XXX3.13} says that $H_n(K;H^p(X;A))=0$ for all $n<p$
and all closed subsets $A$ of $X$.
Thus, $H_n(L;H^p(X;A))=0$ for all $n<p$
and all closed subsets $A$ of $X$. Applying \ref{XXX3.13} again we get
 $SP(L)\in AE(X)$.
\par If $L$ is simply connected, then 2) implies 4) by \ref{XXX2.10}.
 \end{pf}

\begin{Cor} \label{XXX3.18} If $K$ and $L$ are countable, pointed,
connected CW
complexes, then the following conditions are equivalent:
\par{1.} $\dim_G(K)=\dim_G(L)$ for all Abelian groups $G$.
\par{2.} $\dim_Z(K\wedge M)=\dim_Z(L\wedge M)$ for all
pointed CW complexes $M$.
\end{Cor} 
\begin{pf} In view of  \ref{XXX3.16}, Condition 1)
is equivalent to $SP(K)\sim_{Gr} SP(L)$,
and that is equivalent (see \ref{XXX3.17}) to Condition 2).
\end{pf}

\begin{Rem} \label{extrarem}
\ref{XXX3.18} is dual to the well known
characterization of $\dim_G(X)=\dim_G(Y)$ for all Abelian groups $G$
($X$ and $Y$ are compact). Namely, $\dim_G(X)=\dim_G(Y)$ is equivalent
to $\dim_\Z(X\times T)=\dim_\Z(Y\times T)$ for all compact spaces $T$ (see \cite{K}).
\end{Rem}

\section{Cohomological dimension of compact spaces}

The purpose of this section is to show that, given
a compact space $X$, the class $\{SP(L)\mid \ExD(X)\leq SP(L)\}$
has a minimum which we call {\bf the cohomological
dimension} of $X$. To do so we will need
the basics of Bockstein Theory (see \cite{B} or \cite{K}).

\begin{Def} [Bockstein groups] \label{XXX6.1.0}
The set $\B$ of {\bf Bockstein groups}  is 
$$\{\Q\}\cup \bigcup\limits_{\p \ prime} 
\{\Z /\p,\Z /\p^{\infty}, \Z _{(\p)} \},$$
where $\Z /\p^{\infty}$ is the $\p$-torsion of $\Q /\Z$,
and $\Z _{(\p)}$ are the rationals whose
denominator is not divisible by $\p$.
\end{Def}

\begin{Def} [Bockstein basis] \label{XXX6.1.00}
Given an abelian group $G$
its   {\bf Bockstein basis}  $\sigma(G)$ is the subset
of $\B$ defined as follows:
\par{1.} $\Q\in \sigma(G)$ iff $\Q\otimes G\ne 0$, 
\par{2.} $\Z/\p\in \sigma(G)$ iff $(\Z/\p)\otimes G\ne 0$, 
\par{3.} $\Z_{(\p)}\in \sigma(G)$ iff $(\Z/\p^{\infty})\otimes G\ne 0$, 
\par{4.} $\Z/\p^{\infty}\in \sigma(G)$ iff 
$(\Z/\p^{\infty})*G\ne 0$ (here $H*G$ is the torsion product
of groups $H$ and $G$) or $(\Z/\p)\otimes G\ne 0$. 
\end{Def}

\begin{Rem} \label{XXX6.1.00rem}
Our definition of Bockstein basis
is slightly different from the standard ones
(see \cite{K} or  \cite{D$_8$}). Namely, if
 $\Z_{(\p)}\in \sigma(G)$ (respectively,
($\Z/\p\in \sigma(G)$), then $\Z/\p\in \sigma(G)$
(respectively, $\Z/\p^{\infty}\in \sigma(G)$. 
In traditional definitions of Bockstein basis
only one group among $\{\Z_{(\p)},\Z/\p,\Z/\p^{\infty}\}$
is admitted for any $\p$.
Since our only application of Bockstein basis is \ref{XXX6.1.000},
the change of the definition will not cause any problems
as $\dim_{\Z_{(\p)}}(X)\ge \dim_{\Z/\p}(X) \ge \dim_{\Z/\p^{\infty}}(X)$
for all primes $\p$ (see \cite{K}).
\end{Rem}

\begin{Thm} [First Bockstein Theorem] \label{XXX6.1.000}
If $X$ is compact,
then \par\noindent
\centerline {$\dim_G(X)=\max\{\dim_H(X) \mid H\in \sigma(G)\}.$}
\end{Thm}

Here is the existence of cohomological dimension:
\begin{Thm} \label{XXX6.1.0000}
If $X$ is compact,
then there is a countable CW complex $K_X$
such that $SP(K_X)\in AE(X)$
and $SP(K_X)\leq SP(L)$
for every CW complex $L$ satisfying 
$SP(L)\in AE(X)$.
\end{Thm}
\begin{pf} Consider the set $\cal B_X$
of all Bockstein groups $H$ such that $\dim_H(X)<\infty$.
Put $K_X=\bigvee\limits_{H\in\cal B_X}K(H,\dim_H(X))$.
Notice that $K_X\in AE(X)$, hence $SP(K_X)\in AE(X)$
(see \ref{XXX2.10a}). Suppose $SP(L)\in AE(X)$
for some CW complex $L$ and $SP(K_X)\in AE(Y)$
for some compact space $Y$. We need to show
that $SP(L)\in AE(Y)$.
\par Since $SP(K_X)\in AE(Y)$ is equivalent to
$H_i(K_X;H^n(Y,A))=0$ for all $i<n$ and all closed subsets $A$ of $Y$
(see  \ref{XXX3.13}), one has
$H_i(K(H,\dim_H(X));H^n(Y,A))=0$
for all $i<n$ and all closed subsets $A$ of $Y$.
This, in turn, implies $K(H,\dim_H(X))\in AE(Y)$
(see \ref{XXX3.16} and \ref{XXX3.13}) for all $H\in \cal B_X$.
Thus $\dim_H(Y)\leq \dim_H(X)$ for all Bockstein groups $H$
and, in view of \ref{XXX6.1.000}, 
$\dim_G(Y)\leq \dim_G(X)$ for all Abelian groups $G$.
Since, for any compact space $T$, $SP(L)\in AE(T)$
is equivalent to $\dim_{H_i(L)}(T)\leq i$
for all $i$, one has
$\dim_{H_i(L)}(Y)\leq \dim_{H_i(L)}(X)\leq i$ for all $i$,
i.e. $SP(L)\in AE(Y)$.
\end{pf}

\section{Extension types of Eilenberg-MacLane spaces }

The following problem was posed in \cite{D-D$_1$}.


\begin{Problem} \label{XXX6.1} Find necessary and sufficient conditions
for two CW complexes $K$ and $L$ to be of the same extension type.

\end{Problem} 

In this section we consider \ref{XXX6.1} in case of $K$
being the infinite symmetric product and $L$ being an Eilenberg-MacLane space.
 Our characterization of infinite symmetric products having extension type of
an Eilenberg-MacLane space
involves
an enlargement of Bockstein basis.

\begin{Def} \label{XXX6.tauG}
Let $G$ be an abelian group.  $\tau(G)$ is the subset of Bockstein groups
containing $\sigma(G)$ and satisfying the following conditions
for each prime $\p$:
\par{1.} $\Z/\p\in \tau(G)$ iff $\Z/\p^{\infty}\in \sigma(G)$, 
\par{2.} $\Z_{(\p)}\in \tau(G)$ iff $\Z/\p^{\infty}\in \sigma(G)$
and $\Q\in \sigma(G)$. 
\end{Def}

\begin{Lem} \label{XXX6.tauG.1} Suppose $G$ and $F$ are non-trivial Abelian groups
and $m\ge 1$. If $\sigma(F)\setminus \tau(G)\ne\emptyset$,
then there is a compact space $X$ such that $\dim_G(X)=m$
and $\dim_F(X)=\infty$.
\end{Lem} 
\begin{pf} Define $\alpha:\B\to N$ as follows:
$\alpha(H)=\infty$ iff $H\in \B\setminus\tau(G)$,
$\alpha (H)=m+1$ iff $H\in \tau(G)\setminus \sigma(G)$,
and
$\alpha (H)=m$ iff $H\in \sigma(G)$,
\par Let us show that $\alpha$ is a Bockstein function,
i.e. the following inequalities hold for all primes $\p$:
\par 1. $\alpha(\Z/\p^\infty)\leq \alpha(\Z/\p)\leq \alpha(\Z/\p^\infty)+1$,
\par 2. $\alpha(\Z/\p)\leq \alpha(\Z_{(\p)})$,
\par 3. $\alpha(\Q)\leq \alpha(\Z_{(\p)})$,
\par 4. $\alpha(\Z_{(\p)})\leq \max(\alpha(\Q),\alpha(\Z/\p^\infty)+1)$,
\par 5. $\alpha(\Z/\p^\infty)\leq \max(\alpha(\Q),\alpha(\Z_{(\p)})-1)$.
\par Inequalities 1) can fail only if 
$\Z/\p\in\tau(G)$ and $\Z/\p^\infty\in\sigma(G)$. 
In that case $\alpha(\Z/\p^\infty)=m$ and $m\leq \alpha(\Z/\p)\leq m+1$, so Inequalities
1) hold.
Inequality 2) can fail only if 
$\Z_{(\p)}\in\tau(G)$. In that case, however, 
$\Z/\p^\infty\in\sigma(G)$ implying $\Z/\p\in \tau(G)$.
Since $\Z_{(\p)}\in\sigma(G)$ implies $\Z/\p\in \sigma(G)$,
Inequality 2) holds.
Inequality 3) can fail only if 
$\Z_{(\p)}\in\tau(G)$. In that case, however, 
$\Q\in\sigma(G)$, so
Inequality 3) holds.
Inequality 4) can fail only if $\Z/\p^\infty\in\tau(G)$
and $\Q\in\tau(G)$. That however implies 
$\Z/\p^\infty\in\sigma(G)$
and $\Q\in\sigma(G)$. Consequently, $\Z_{(\p)}\in\tau(G)$ and
$\alpha(\Z_{(\p)})\leq m+1\leq \max(\alpha(\Q),\alpha(\Z/\p^\infty)+1)$,
i.e. 4) holds.
Inequality 5) can fail only if $\Z_{(\p)}\in \tau(G)$.
Therefore $\Z/\p^\infty\in\sigma(G)$ and $\Q\in\sigma(G)$.
Hence $\alpha(\Z/\p^\infty)=m\leq \max(\alpha(\Q),\alpha(\Z_{(\p)})-1)$.
\par By Dranishnikov's Realization Theorem (see \cite{D$_1$} or \cite{D$_8$})
there is a compactum $X$
such that $\dim_H(X)=\alpha (H)$ for all $H\in \B$.
It is clear, in view of Bockstein's First Theorem \ref{XXX6.1.000}, that $X$
satisfies the desired conditions.
\end{pf}

\begin{Lem} \label{XXX6.tauG.2} Suppose $G$ and $F$ are non-trivial Abelian groups
and $m\ge 1$. If $\sigma(F)\subset \tau(G)$ and
$\sigma(F)\setminus \sigma(G)\ne\emptyset$,
then there is a compact space $X$ such that $\dim_F(X)=\dim(X)=m+1$
and $\dim_G(X)=m$.
\end{Lem} 
\begin{pf} Notice that $\sigma(F)\setminus \sigma(G)$
contains only groups of the form $\Z/\p$ or $\Z_{(\p)}$
for some $\p$.
\par
Case I: There is a prime $\q$ such that
$\Z/\q\in \sigma(F)\setminus \sigma(G)$.
Define $\alpha:\B\to N$ 
by sending $\Z/\q$ and $\Z_{(\q)}$ to $m+1$, and sending
all the other groups to $m$.
 Notice that $\alpha$ is a Bockstein function.
  By Dranishnikov's Realization Theorem (see \cite{D$_1$} or \cite{D$_8$})
there is a compactum $X$
such that $\dim_H(X)=\alpha (H)$ for all $H\in \B$
and $\dim(X)=\max(\alpha )$.
It is clear, in view of Bockstein's First Theorem, that $X$
satisfies the desired conditions.
\par Case II: There is a prime $\q$ such that
$\Z_{(\q)}\in \sigma(F)\setminus \sigma(G)$
and all groups in $\sigma(F)\setminus \sigma(G)$
are torsion-free.
Define $\alpha:\B\to N$ 
by sending $\Z_{(\q)}$ to $m+1$, and sending
all the other groups to $m$. Notice that $\alpha$ is a Bockstein function.
 By Dranishnikov's Realization Theorem (see \cite{D$_1$} or \cite{D$_8$})
there is a compactum $X$
such that $\dim_H(X)=\alpha (H)$ for all $H\in \B$
and $\dim(X)=\max(\alpha )$.
It is clear, in view of Bockstein's First Theorem, that $X$
satisfies the desired conditions.
\end{pf}

\begin{Thm} \label{XXX6.1.1} Suppose $L$ is a pointed countable
CW complex, $G$ 
is an Abelian group, and $n\ge 1$.
 $SP(L)$ is of the same extension type as $K(G,n)$
if and only if the following conditions are satisfied:
\par{a.} $H_i(L)=0$ for $i<n$.
\par{b.} $\sigma(H_n(L))=\sigma(G)$.
\par{c.} $\sigma(H_i(L))\subset \tau(G)$ for all $i\ge n$.
\end{Thm} 
\begin{pf} We can reduce the general case to that of $K(G,n)$
being a countable CW complex. Indeed,  \ref{XXX6.1.000}
implies that any $K(G,n)$ has the extension type
of $K(G',n)$ such that $G'$ is countable and $\sigma(G)=\sigma(G')$.
\par
Assume $SP(L)$ is of the same extension type as $K(G,n)$.
Since $L\sim_{Gr} K(G,n)$ (see  \ref{XXX3.16} and  \ref{XXX3.17}), a) follows.
Pick $i\ge n$ and denote $H_i(L)$ by $F$.
Suppose $\sigma(F)\setminus \tau(G)\ne \emptyset$.
By \ref{XXX6.tauG.1} there is a compactum $X$ such that $\dim_G(X)=n$
and $\dim_F(X)=\infty$. Since $SP(L)$ is of the same extension type as $K(G,n)$,
$\dim_G(X)=n$ implies $SP(L)\in AE(X)$. Consequently,
$\dim_F(X)=\dim_{H_i(L)}(X)\leq i$ (see  \ref{XXX2.10b}), a contradiction.
Thus c) holds.
\par Denote $H_n(L)$ by $F$. Thus $\sigma(F)\subset \tau(G)$.
Suppose $\sigma(F)\setminus \sigma(G)\ne \emptyset$.
By \ref{XXX6.tauG.2} there is a compact space $X$ such that 
$\dim_G(X)=\dim(X)=n+1$
and $\dim_F(X)=n$.
Since $\dim(X)=n+1$, $\dim_{H_i(L)}(X)\leq i$ for all $i\ge n$.
Consequently, see  \ref{XXX2.10b}, $SP(L)\in AE(X)$
which implies $K(G,n)\in AE(X)$ contradicting $\dim_G(X)=n+1$.
That proves $\sigma(F)\subset \sigma(G)$.

\par Suppose $\sigma(G)\setminus\sigma(F)\ne\emptyset$.
By \ref{XXX6.tauG.2} there is a compact space $X$ such that $\dim_F(X)=\dim(X)=n+1$
and $\dim_G(X)=n$. Since $SP(L)$ is of the same extension type as $K(G,n)$,
$\dim_G(X)=n$ implies $SP(L)\in AE(X)$. Consequently,
$\dim_F(X)=\dim_{H_n(L)}(X)\leq n$ (see  \ref{XXX2.10b}), a contradiction.
That proves b).
\par Suppose a), b), and c) hold.
 If $K(G,n)\in AE(X)$ (i.e., $\dim_G(X)\leq n$),
then $\dim_F(X)\leq n+1$ for all $F\in\tau(G)$ in view of Bockstein's Inequalities.
Hence, $\dim_{H_i(L)}(X)\leq i$ for all $i\ge n$
which implies $SP(L)\in AE(X)$ (see  \ref{XXX2.10b}). 
That shows $K(G,n)\leq SP(L)$.
\par If $SP(L)\in AE(X)$, then $K(H_n(L),n)\in AE(X)$ (see  \ref{XXX2.10b})
which is equivalent to $K(G,n)\in AE(X)$ in view of
$\sigma(G)=\sigma(H_n(L))$ and the First Bockstein Theorem.
\end{pf}

\section{Dimension types and the connectivity index}

\par Shchepin's connectivity index is of homological nature.
Similarly, one can introduce the {\bf homotopy connectivity index}
$hcin(K)$.

\begin{Def} \label{XXXaa3.1} Suppose $X$ is a pointed
 space. $hcin(X)$ is a non-negative integer defined as follows:
\par a. $hcin(X)=0$ means that $X$ is not path-connected.
\par b. If $\infty > r > 0$, then
$hcin(X)=r$ means that $\pi_r(X)\ne 0$
and $\pi_k(X)=0$ for all $0\leq k<r$. 
\par c. $hcin(X)=\infty$ means that $\pi_k(X)=0$ for all $0\leq k< \infty$.
\end{Def} 

Homotopy connectivity index can be easily dualized. Following
G.Whitehead \cite{Wh} (p.421--423) we introduce
the {\bf anticonnectivity index} $acin(X)$ as follows.

\begin{Def} \label{XXXaaa3.1} Suppose $X$ is a pointed
 space. $acin(X)$ is an integer greater than or equal to $-1$
or infinity defined as follows:
\par a. $acin(X)=-1$ means that $X$ is path-connected and all its
homotopy groups are trivial.
\par b.
If $0 \leq r <\infty$, then
$acin(X)=r$ means that $\pi_r(X)\ne 0$
and $\pi_k(X)=0$ for all $k>r$. 
\par c.
$acin(X)=\infty$ means that infinitely many homotopy groups of $X$
are non-trivial.
\end{Def}

\par We start with the the concept of the total function
space which is related to Shchepin's \cite{S} concept of the total
cohomology of a space.

\begin{Def} \label{XXXa3.1} Suppose $X$ is a pointed
compact space and $P$
is a pointed CW complex. The {\bf total function space} $Tot(P^X)$
is the wedge of all function spaces $P^{X/A}$, where $A$ 
is a closed subspace of $X$.
\end{Def}

Using homotopy connectivity index and the concept of the total
function space one can introduce a new relation on the class
of compact spaces: $X\sim Y$ iff $hcin(Tot(P^X))=hcin(Tot(P^Y))$
for all pointed CW complexes $P$. It turns out this relation
means exactly that $\Sigma X$ and $\Sigma Y$ are of the same
extension dimension. The first part of this section is devoted
to that fact and culminates in \ref{XXX.1}.
\par In the case of pointed CW complexes one can define
a relation $K\sim L$ to mean that $cin(K\wedge M)=cin(L\wedge M)$
for all $M$. We saw in \ref{XXX3.17} and \ref{XXX3.18}
that, in the case of countable CW complexes, that relation is identical
with equality of extension types of $SP(K)$ and $SP(L)$.
Dualizing that relation leads to $K\sim L$ iff
$acin(P^K)=acin(P^L)$ for all $P$. In the last part of this section
we investigate how extension types are connected to that relation.

\begin{Lem} \label{XXXhcin.1} Let $X$ be a compact connected space.
If $P$ is a pointed CW complex, then
$$hcin(Tot(\Omega P)^X)\ge hcin(Tot( P)^X)-1.$$
Moreover, if $hcin(Tot( P)^X)\ge 1$, then
$$hcin(Tot(\Omega P)^X)=hcin(Tot( P)^X)-1.$$
\end{Lem} 
\begin{pf} If $hcin(Tot( P)^X)=0$, then the inequality obviously
holds, so assume $hcin(Tot( P)^X)$ is at least $1$. Now, 
\ref{XXXhcin.1} follows from equality
$[S^{n-1},(\Omega P)^{X/A}]=[S^n,P^{X/A}]$
for all closed subsets $A\ne\emptyset$ of $X$ and all $n\ge 1$.
\end{pf}

\begin{Rem} \label{XXXhcin.2} If $hcin(Tot( P)^X)=0$,
then $hcin(Tot(\Omega P)^X)$ may be arbitrarily high.
For example, if $P=X=S^1$, then $Tot(\Omega P)^X$ is contractible
but $Tot( P)^X$ is not path-connected.
\end{Rem} 

\begin{Lem} \label{XXXhcin.3} Let $X$ be a compact connected space and $k > 0$.
If $P$ is a pointed CW complex, then the following conditions are equivalent:
\par a. $hcin(Tot( P)^X)\ge k$.
\par b. $\Omega^k (P)\in AE(X)$.
\end{Lem} 
\begin{pf} Notice that both $hcin(Tot( P)^X)$
and $\Omega^k (P)$ depend only on the component of the base point
of $P$, so we may reduce \ref{XXXhcin.3}
to the case of $P$ being connected.
\par
Special Case: $k=1$. Notice that
$hcin(Tot(P^X))\ge 1$ is equivalent to
$[X/A,P]=\ast$ for all closed subsets $A\ne\emptyset$ of $X$.
That statement, in view of \ref{XXX2.11}, is equivalent
to $\Omega P\in AE(X)$.
\par If $k > 1$, then $hcin(Tot(P^X))\ge k$ is equivalent
(in view of  \ref{XXXhcin.1})
to $$hcin(Tot((\Omega^{k-1}P)^X))\ge 1$$ which is equivalent,
by the Special Case, to $\Omega^k(P)\in AE(X)$.
\end{pf}

\begin{Thm} \label{XXX.1} If $X$ and $Y$ are non-empty connected
compact spaces, the following
conditions
are equivalent:
\par a. $\Sigma X$ and $\Sigma Y$ are of the same
extension dimension.
\par{b.} $\Omega(P)\in AE(X)$ is equivalent to
 $\Omega(P)\in AE(Y)$ for all pointed CW complexes $P$.
\par{c.} $hcin(Tot(P^X))=hcin(Tot(P^Y))$
for all pointed CW complexes $P$.
\end{Thm}

\begin{pf}  a)$\implies$ b). 
Suppose $\Omega P\in AE(X)$ for some pointed CW complex $P$.
Since $\Omega P=\Omega P_0$, where $P_0$ is the component
of the base point of $P$,
we may assume $P$ is connected.
By \ref{XXX2.11}, $P\in AE(\Sigma X)$ which implies
$P\in AE(\Sigma Y$, as $\Sigma Y$ is of the same extension
dimension as $\Sigma X$.
Again, by \ref{XXX2.11}, $\Omega P\in AE(X)$.
The same argument shows that  $\Omega P\in AE(Y)$
implies  $\Omega P\in AE(X)$ for any pointed CW complex $P$.

\par
b)$\implies$ c). Suppose $hcin(Tot(P^X))\ge k$.
It suffices to show (by symmetry) that
$hcin(Tot(P^Y))\ge k$.
If $k=0$ it is obviously true, so assume $k\ge 1$.
However, in that case,
$hcin(Tot(P^T))\ge k$ is equivalent (see \ref{XXXhcin.3})
to $\Omega^k(P)\in AE(T)$ for any compact connected space $T$,
so b)$\implies $ c) follows from \ref{XXXhcin.3}
as $\Omega^k(P)\in AE(X)$ is equivalent to $\Omega^k(P)\in AE(Y)$.
\par c)$\implies$ a).
Suppose $P\in AE(\Sigma X)$. It suffices to show $P\in AE(\Sigma Y)$.
Since $\Sigma X$ is connected and contains at least two points,
$P$ must be connected (otherwise we pick $x_1\ne x_2\in \Sigma X$,
map them to two different components of $P$, and that map
cannot be extended over $\Sigma X$), so (see \ref{XXX2.11})
$P\in AE(\Sigma X)$ is equivalent to $\Omega P\in AE(X)$.
Hence $\Omega P\in AE(X)$ and, by \ref{XXX2.11},
$P\in AE(\Sigma Y)$.
 \end{pf}

 \ref{XXX.1} implies that if two compacta $X$ and $Y$ are of the same
dimension type, then the total function spaces $Tot(P^X)$ and $Tot(P^Y)$
have the same homotopy connectivity index. The rest of this section
is devoted to the dual result: if two countable complexes $K$ and
$L$ are of the same extension type, then the function spaces
$P^K$ and $P^L$ have the same anticonnectivity index for certain $P$.

\par Since some of the techniques of this section are well-known
(see \cite{D-W$_3$} and \cite{D$_8$}), we will only outline
how one translates known results in terms of trunctated cohomology
to results in terms of function spaces.
\par The following result corresponds to the Combinatorial Vietoris-Begle
Theorem of \cite{D-W$_3$} (see Lemma 2 there or Lemma 5.9 in \cite{D$_8$})
and the proof is similar.


\begin{Prop} \label{XXX.2}
Suppose $P$ is a CW complex and $f:K\to L$ is a combinatorial map from a CW
complex $K$
to a finite simplicial complex $L$.
If $P^{f\1(\Delta)}$ is weakly contractible for each simplex $\Delta$
of $L$, then $f^*:P^L\to P^K$ is a weak homotopy equivalence.

\end{Prop}

The next two results extract essential parts of Lemma 3 in of \cite{D-W$_3$}
with proofs being similar.


\begin{Lem} \label{XXX.3}
Suppose $P$ is a CW complex. The following conditions are equivalent:
\par{1.} $[S^k,P]$ is finite for each $k\ge 0$ (unpointed spheres).
\par{2.} $[K,P]$ is finite for each finite CW complex $K$.
\par{3.} for each pair $(K,L)$ of finite CW complexes
and each map $f:L\to P$, the set
of all homotopy classes rel.$L$ of maps $g:K\to P$ such that $g|L=f$, is
finite.

\end{Lem}


\begin{Lem} \label{XXX.4}
Suppose $P$ is a CW complex such that
 $[S^k,P]$ is finite for each $k\ge 0$ (unpointed spheres).
If $K$ is a countable CW complex and $f:K\to P$ is a map such that $f|L\sim
0$ for each finite subcomplex
$L$ of $K$, then $f\sim const$.

\end{Lem}


\begin{Thm} \label{XXX.5}
Suppose $P$ is a CW complex and $K$ is a countable CW complex
such that $P^K$ is weakly contractible.
Suppose $L$ is a countable complex and $f:L\to P$ is a homotopically
non-trivial map.
There is a compactum $X$ and a map $g:X\to L$ so that $f\circ g$ is 
homotopically non-trivial
and $K\in AE(X)$ if one of the following conditions is satisfied:
\par{1.} $K$ is compact and $L$ is finitely dominated,
\par{2.} $[S^k,P]$ is finite for each $k\ge 0$ (unpointed spheres).

\end{Thm}

\begin{pf}  In case of 2), one follows the same technique as described
in \cite{D$_8$} (see 5.5-5.7 there). We will outline how to use that technique
in case of both $K$ and $L$ being compact.
The inductive step consists of a homotopically non-trivial map $f':L'\to P$
and a map $h':L_0\to K$, where $L'$ is a compact simplicial complex
and $L_0$ is a subcomplex of $L'$. As in Lemma 5.6 of \cite{D$_8$}
one constructs a combinatorial map $\pi:L''\to L'$ 
so that $\pi\1(\Delta)$ is either contractible or homotopy
equivalent to $K$
for each simplex $\Delta$ of $L'$. Moreover, the composition
$\pi\1(L_0)\to L_0\to K$ extends over $L''$. 
A construction in \cite{Dy$_4$} consists of $L''$ being a subcomplex
of $L'\times K$, \cite{D$_8$} constructs $L''$ as the pull-back of
a certain diagram. \ref{XXX.2} says that the composition $L''\to L'\to P$
is homotopically non-trivial. Starting with $f:L\to P$ one can
construct inductively finite simplicial complexes $L_n$
and maps $f_n:L_{n+1}\to L_n$ so that the inverse limit $X$
of the sequence $\ldots L_{n+1}\to L_n\ldots $ satisfies
$K\in AE(X)$, $L_1$ is homotopically equivalent to $L$,
and the composition $L_n\to\ldots\L_1\to P$ is homotopically
non-trivial for all $n$ (see 5.5 and 5.9 of \cite{D$_8$}
for details).
\par If $L$ is finitely dominated, we pick $L'$ finite
and maps $u:L\to L'$, $d:L'\to L$ so that $d\circ u\approx id_L$.
Let $f'=f\circ d:L'\to P$. Notice that $f'$ is homotopically non-trivial.
By the previous case, there is compactum $X$ and a map $g':X\to L'$
so that $K\in AE(X)$ and $f'\circ g'$ is homotopically non-trivial.
Let $g=d\circ g'$. Since $f\circ g=f'\circ g'$, it is homotopically non-trivial.
 \end{pf}


\begin{Cor} \label{XXX.6} Suppose $K\leq L$ are countable pointed
complexes.
Suppose $P$ is a pointed CW complex so that $\pi_i(P^K)=0$ for all $i\ge r$,
where $r\ge 1$. Then the homotopy groups $\pi_i(P^L)$ are trivial for all $i\ge r$ if one
of the following conditions hold:
\par{1.} $\pi_i(P)$ is finite for all $i\ge r$,
\par{2.} $K$ is compact and $L$ is finitely dominated,
\par{3.} both $K$ and $L$ are finitely dominated and $r\ge 2$.

\end{Cor}

\begin{pf}  Suppose $\pi_m(P^L)\ne 0$ for some $m\ge r$.
Put $L'=\Sigma(L)$ in cases 1) and 2), put $L'=\Sigma^2(L)$
in case 3).
Put $P'=\Omega^{m-1}(P)$ in cases 1) and 2), put $P'=\Omega^{m-2}(P)$
in case 3). Put $K'=K$ in cases 1) and 2), put $K'=\Sigma(K)$
in case 3).
Notice that $K'$ and $L'$ are homotopically equivalent to
compact CW complexes in cases 2) and 3) (see \cite{Wall}).
In all cases there is a homotopically non-trivial map $f:L'\to P'$,
so there is a map $g:X\to L'$ such that $K'\in AE(X)$
and $f\circ g$ is homotopically non-trivial (see \ref{XXX.5}). 
This contradicts \ref{XXX2.1}. Indeed, $K\leq L$ implies $\Sigma(K)\leq \Sigma(L)$
(see \cite{D-D$_1$}). That means $K'\in AE(X)$ implies that
any map $X\to L'$ is null-homotopic by \ref{XXX2.1}.
 \end{pf}

\begin{Cor} \label{XXX.6b} Suppose $K\leq L$ are countable pointed
complexes and $\p$ is a prime.
If $H^\ast(K;\Z/\p)=0$, then $H^\ast(L;\Z/\p)=0$.
\end{Cor}
\begin{pf} Suppose $H^\ast(K;\Z/\p)=0$ and $H^n(L;\Z/\p)\ne 0$
for some $n\ge 0$. If $n=0$, then $L$ must be disconnected
and $K$ is connected. Hence $K\in AE(S^1)$ and $L\notin AE(S^1)$,
a contradiction. Thus, $n\ge 1$.
Put $P=K(\Z/\p,n+1)$. Notice that $\pi_1(P^L)=H^n(L;\Z/\p)\ne 0$
and $\pi_i(P^K)=H^{n+1-i}(K;\Z/\p)=0$ for all $i\ge 1$,
a contradiction in view of \ref{XXX.6}.
\end{pf}


\begin{Thm} \label{XXX.7} Suppose $L$ is a countable, finite-dimensional CW
complex
and $\p$ is a prime number. If $H^n(L;\Z/\p)\ne 0$ for some $n\ge 1$,
then there is a non-trivial map $f:X\to\Sigma^{n+2}(L)$ from a compactum $X$
so that $\dim_{\Z[1/\p]}X=1=\dim_{\Z/\p}X$.

\end{Thm}
\begin{pf}  Let $M=M(\Z/\p,1)$ and let $L'=\Sigma^{n+2}(L)$.
In particular, $H^{2n+2}(L';\Z/\p)\ne 0$. Notice that $K(\Z,2n+3)^M$
is a $K(\Z/\p,2n+2)$ so there is a non-trivial map
$g:L'\to K(\Z,2n+3)^M$. Its adjoint $g':L'\wedge M\to K(\Z,2n+3)$
is non-trivial and we may assume that its image is contained
in a finite subcomplex $A$ of $K(\Z,2n+3)$ which is $(2n+2)$-connected.
Notice that $P=A^M$ is simply connected and its homotopy groups are
finite. Indeed, the homotopy groups of $P$ coincide with groups of
homotopy clases of maps from suspensions of $M$ to $A$.
Those sets are the same as mod $\p$ homotopy groups of $A$
(see \cite{N}). In view of Proposition 1.4 on p.3 in \cite{N} 
one has an exact sequence
$$ 0\to \pi_m(A)\otimes \Z/\p\to \pi_m(A;\Z/\p)\to \pi_{m-1}(A)\ast \Z/\p\to 0$$
and, since the homotopy groups of $A$ are finitely generated,
the homotopy groups of $P$ are finite.
By Miller's Theorem (Sullivan Conjecture - see \cite{Mi}),
$A^{K(\Z/\p,1)}$ is weakly contractible as $A$ is finite-dimensional. Therefore,
$P^{K(\Z/\p,1)}=(A^{K(\Z/\p,1)})^M$ is weakly contractible.
Since
$K(\Z[1/\p],1)\wedge M$ is contractible (compute its homology groups),
we get that $P^{K(\Z[1/\p],1)}$ is weakly contractible.
Let $K=K(\Z/\p,1)\vee K(\Z[1/\p],1)$. Notice that $P^K$
is weakly contractible. Applying \ref{XXX.5} one gets a non-trivial map
$f:X\to L'$ so that $K\in AE(X)$.
 \end{pf}

\ref{XXX.7} and \ref{XXX2.1} imply the following.


\begin{Cor} \label{XXX.8} Suppose $L$ is a countable, finite-dimensional CW
complex
and $\p$ is a prime number. If $H^n(L;\Z/\p)\ne 0$ for some $n\ge 1$,
then there is a compactum $X$
so that $\dim_{\Z[1/\p]}X=1=\dim_{\Z/\p}X$ and $L$ is not an absolute
extensor of $X$.

\end{Cor}

\section{Extension types of infinite symmetric products }

In this section we consider Problem \ref{XXX6.1}
in the case of $K$ being an infinite symmetric space
and $L$ being a compact CW complex (respectively, a countable,
finite-dimensional CW complex).

\begin{Lem} \label{XXX6.1.5} Suppose $\p$ is a prime. The following conditions are
equivalent
for any countable connected pointed CW complex $K$:
\par{1.} $H^*(SP(K);\Z/\p)=0$.
\par{2.} $\Z/\p^\infty\notin\sigma(H_s(K))$ for all $s\ge 1$.

\end{Lem}
\begin{pf} 1)$\implies$ 2). Let $G_i=H_i(K)$ for $i\ge 1$.
If $\p\cdot G_i\ne G_i$ for some $i$,
there is a non-trivial map $K(G_i,i)\to K(G_i/\p\cdot G_i,i)$,
a contradiction as $G_i/\p\cdot G_i$ is a direct sum of copies of $\Z/\p$.
Thus $\p\cdot G_i=G_i$ for all $i$.
If $\p-\Tor(G_i)\ne 0$ for some $i$,
then $K(G_i,i)$ dominates $K(\Z/\p^\infty,i)$. To complete the proof of 
1)$\implies$ 2) it suffices to show that
$K=K(\Z/\p^\infty,i)$ has non-trivial
$\Z/\p$-cohomology. Put $L=K(\Z/\p,i+2)$
an notice that $K\leq L$ in view of Bockstein's Inequalities. 
Put $P=K(\Z/\p,i+4)$
and notice that the triviality of $\Z/\p$-cohomology of $K$
means $\pi_n(P^K)=0$ for all $n\ge 1$. \ref{XXX.6b} says that
$\pi_n(P^L)=0$ for all $n\ge 1$, a contradiction as $\pi_2(P^K)=\Z/\p$.
\par 2)$\implies$ 1). Notice that $M(\Z[1/\p],1)\leq_{Gr} K$.
Indeed, $\Z/\p^\infty\notin\sigma(H_s(K))$ for all $s\ge 1$ means
that $H_*(K)$ has no $\p$-torsion and is divisible
by $\p$. If $F\otimes \Z[1/\p]=0$ and $F$ is a Bockstein group, then 
$F$ must be either $\Z/\p$ or $\Z/\p^\infty$. Therefore
$H_*(K;F)=0$ which completes the proof of $M(\Z[1/\p],1)\leq_{Gr} K$.
Also, $H^*(K(\Z[1/\p],1);\Z/\p)=0$. If $H^*(SP(K);\Z/\p)\ne 0$,
then there is a compactum $X$ so that $K(\Z[1/\p],1)\in AE(X)$
but $SP(K)\notin AE(X)$ (see \ref{XXX.8})
which contradicts
\par\noindent $K(\Z[1/\p],1)\sim_{Gr} M(\Z[1/\p],1)\leq_{Gr}
K\sim_{Gr}SP(K)$
(see \ref{XXX3.16} and \ref{XXX3.17}).
 \end{pf}


\begin{Thm} \label{XXX6.2} Suppose $L$ is a connected countable
CW complex and $n\ge 0$.
If $\Sigma^n(SP(L))$ has the extension type of a countable,
finite-dimensional and non-trivial CW complex,
then either $SP(L)$ is of the same extension type as $K(\Z_{(l)},1)$
for some subset $l$ of primes or 
it is of the same extension type as
$K(\Q,m)$ for some $m\ge 1$.
\end{Thm}
\begin{pf} 
Let $G_i=H_i(L)$ for $i\ge 1$.
\par Suppose $SP(L)$ is not of the same extension type as
$K(\Q,m)$ for all $m$.
There is the smallest $r$ with $\sigma(G_r)\ne\sigma(\Q)$ (see \ref{XXX6.1.1}).
There must be a prime $\p$ with $\Z/\p^\infty\in\sigma(G_r)$ which means
$H^*(SP(L);\Z/\p)\ne 0$ (see \ref{XXX6.1.5}).
Suppose $\Sigma^n(SP(L))\sim K$, $K$ is finite-dimensional.
That means $H^*(K;\Z/\p)\ne 0$ (see \ref{XXX.6b}) and there is a
compactum $X_\p$ with $\dim_{\Z[1/\p]}(X_\p)=1=\dim_{\Z/\p}(X_\p)$ so that
$K$ is not an
absolute extensor of $X_\p$ (see \ref{XXX.8}).
If $r\ge 2$, then $\dim_{G_i}(X_\p)\leq i$ for all $i$ as $\dim_\Z(X_\p)\leq 2$.
Hence $SP(L)\in AE(X_\p)$ which implies $K\in AE(X_\p)$, a
contradiction.
Thus, $r=1$. If $\sigma(G_1)\ne \sigma(\Z_{(l)})$ for all sets of primes $l$,
then $\p$ above may be chosen so that $\Z_{(p)}\notin\sigma(G_1)$
which implies $\dim_{G_1}(X_\p)=1$. Again, $SP(L)\in AE(X_\p)$ which
implies $K\in AE(X_\p)$ (see \ref{XXX3.4}), a contradiction.
 Assume $\sigma(G_1)=\sigma(\Z_{(l)})$. If $\sigma(G_s)\subset\sigma(\Z_{(l)})$
for each $s>1$, then we are done by  \ref{XXX6.1.1}.
Suppose $\sigma(G_s)\subset\sigma(\Z_{(l)})$ does not hold
for some $s>1$.
There must be a prime $\p\notin l$ so that $\Z/\p^\infty\in\sigma(G_s)$.
Again, there is a compactum $X_\p$ with
$\dim_{\Z[1/\p]}(X_\p)=1=\dim_{\Z/\p}(X_\p)$ so
that $K$ is not an
absolute extensor of $X_\p$.
This implies $\dim_{G_i}X_\p\leq 1$ for all $i$ and $SP(L)\in AE(X_p)$.
Again, $K\in AE(X_\p)$, a contradiction.
 \end{pf}


\begin{Thm} \label{XXX6.3} Suppose $L$ is a connected countable
CW complex and $n\ge 0$.
If $\Sigma^n(SP(L))$ is
 of a compact non-trivial extension type
(i.e., there is a compact CW complex $K$ of the same extension type as
$\Sigma^n(SP(L))$),
then $SP(L)$ is of the same extension type as $S^1$.

\end{Thm}

\begin{pf}  Let $K$ be a compact CW complex of the same extension type as
$\Sigma^n(SP(L))$.
 By \ref{XXX6.2} either
$SP(L)\sim K(\Z_{(l)},1)$ or $SP(L)\sim K(\Q,m)$ for some
$m\ge 1$.
I $l$ is not the set of all primes or $SP(L)\sim K(\Q,m)$ for some
$m\ge 1$,
then there is a prime $\p$ such that $H^*(SP(L);\Z/\p)=0$
(choose $\p\notin l$ if $SP(L)\sim K(\Z_{(l)},1)$ or
any $\p$ if $SP(L)\sim K(\Q,m)$). That implies
$H^*(K;\Z/\p)=0$ by \ref{XXX.6b}. Since $K$ is finite, $\tilde H_*(K)$ must be a torsion
graded group
and $\tilde H_*(K;\Q)=0$. Thus,
$\tilde H_*(SP(L);\Q)=0$ (see \ref{XXX3.17}), a contradiction
as $SP(L)\sim K(\Z_{(l)},1)$ or $SP(L)\sim K(\Q,m)$ for
some $m\ge 1$.
 \end{pf}

Theorem \ref{XXX6.3} generalizes all the known theorems related to the difference
between
extension types:
\par{1.} $S^n$ and $K(\Z,n)$ are of different extension types
for $n\ge 3$ (\cite{D$_2$}).
\par{2.} $S^n$ and $K(\Z,n)$ are of different extension types
for $n\ge 2$ (\cite{D-W$_3$}).
\par{3.} $M(\Z/\p,n)$ and $K(\Z/\p,n)$ are of different extension types
for $n\ge 1$
(\cite{M$_2$}).
\par{3.} $M(\Z/2,1)$ and $K(\Z/2,1)$ are of different extension types
(\cite{L}).
\par{4.} $RP^n$ and $RP^{\infty}$ are of different extension types
for $n\ge 1$
(\cite{D-R$_2$}).

\par Besides generalizing the above-mentioned results, a major reason
for \ref{XXX6.3} was our interest in pursuing ways of proving/disproving
existence of universal spaces of given cohomological dimension.
In \cite{Dy$_3$}, the author generalized a result of Shvedov which states that,
for any compact CW complex $K$, the class of compacta $X$ such that
$K$ is an absolute extensor of $X$ has a universal space.
That generalization deals with $K$ homotopy dominated by a compact
CW complex. Thus a natural way to see if there is a universal space
of a given cohomological dimension is to verify if a particular
CW complex has extension type of a compact CW complex. Theorem \ref{XXX6.3}
closes that route of proving existence of universal spaces
for any infinite symmetric product but $S^1$.


\begin{Lem} \label{XXX6.4} $M(\Q,n)$ and $K(\Q,n)$
are of the same extension type for all $n\ge 2$.

\end{Lem}

\begin{pf}  $M=M(\Q,n)$ can be realized
as the telescope of $f_m:S^n\to S^n$,
where $f_m$ is of degree $m!$ for $m\ge 1$.
In particular, homotopy groups of $M(\Q,n)$ are torsion groups
for $i\ge 2n$.
By Sullivan's Theorem \cite{Su}
there is a map $f:M(\Q,n)\to K$ such that
$H_i(K)=H_i(M)\otimes\Q$, $\pi_i(K)=\pi_i(M)\otimes\Q$
for all $i\ge 1$ and
$f_*:H_i(M)\to H_i(K)$ corresponds to $H_i(M)\otimes\Z\to H_i(M)\otimes \Q$
for all $i\ge 1$. Thus, $f$ is a homotopy equivalence.
In particular, $\pi_i(M)=0$ for $i\ge 2n$.
 Thus, $M(\Q,n)\sim K(\Q,n)$
(see \cite{Dy$_1$}).
 \end{pf}


\begin{Thm} \label{XXX6.5} Suppose $G$ is a countable Abelian group. The
following conditions
are equivalent:
\par{1.} There is a Moore space $M(G,n)$ of the same extension type as
$K(G,n)$.
\par{2.} Either $n=1$ and there is a subset $l$ of primes
such that $K(G,1)$ is of the same extension type as $K(\Z_{(l)},1)$ or $n\ge 2$ and
$K(G,n)$ is of the same extension type as $K(\Q,n)$.
\end{Thm}

\begin{pf}  2)$\implies$ 1) follows from \ref{XXX6.4} and the fact that
one can choose $M(\Q,1)$ to be $K(\Q,1)$.
\par 1)$\implies$ 2). Notice that $\Sigma(M(G,n))$ is homotopy equivalent to
a finite-dimensional CW complex and has the same extension type as
$\Sigma(K(G,n))$. Use \ref{XXX6.2}.
 \end{pf}

\begin{Rem} \label{XXX6.6bef} \ref{XXX6.5} solves the following problem posed by 
Dranishnikov (see \cite{S}, p.983):
\par Is it true that for any compactum $X$ and any countable group $G$
the conditions $M(G,1)\in AE(X)$ and $K(G,1)\in AE(X)$
are equivalent?
 
\end{Rem}


\begin{Problem} \label{XXX6.6} Suppose the extension type of
a countable CW complex $K$ is at most the extension type of
a compact non-contractible CW complex $L$. Is there
a finite-dimensional, countable CW complex $M$ of the same
extension type as $K$?

\end{Problem}

\begin{Rem} \label{XXX6.6b} One cannot replace compactness of $L$
by finite-dimensionality of $L$. Indeed, $K(\Z,2)\leq M(\Q,2)$
but $K(\Z,2)$ does not have the extension type
of a countable finite-dimensional CW complex (see \ref{XXX6.2}
and  \ref{XXX6.1.1}).
 
\end{Rem}


\begin{Problem} \label{XXX6.7} Suppose the extension type of
a countable CW complex $K$ is at most the extension type of
a compact non-contractible CW complex $L$. Is there
a universal space among all compacta $X$ so that $K\in AE(X)$?

\end{Problem}

\medskip
\medskip
\medskip
\medskip
Jerzy Dydak,
Math Dept, University of Tennessee, Knoxville, TN 37996-1300, USA,
E-mail address: dydak@@math.utk.edu

\end{document}